\documentclass[preprint,12pt]{elsarticle}
\usepackage{mathrsfs}
\usepackage{amsfonts}
\usepackage{amssymb,bm}
\usepackage{amsmath}
\usepackage{amsthm}
\usepackage{bbding}
\usepackage{tikz}
\usetikzlibrary{matrix,shapes,snakes}
\allowdisplaybreaks
\newtheorem{thm}{Theorem}[section]

\newtheorem{lem}[thm]{Lemma}
\newtheorem{cor}[thm]{Corollary}

\newtheorem{rmk}[thm]{Remark}

\newproof{pf}{Proof}
\newcommand{\qedd}{\hspace*{\fill}$\Box$\medskip}   

\def\nor{{\rm{N}}}

\journal{}
\begin{document}

\begin{frontmatter}

\title{The compositional inverses of linearized permutation binomials over finite fields}

\author{Baofeng Wu\corref{cor1}}
\ead{wubaofeng@iie.ac.cn}

\cortext[cor1]{Corresponding author. Fax: +86 13426076355.}

\address{State Key Laboratory of Information Security, Institute of Information Engineering, Chinese Academy of Sciences, Beijing 100093, China}

\begin{abstract}
Let $q$ be a prime power and $n$ and $r$ be positive integers. It is
well known that the linearized binomial
$L_r(x)=x^{q^r}+ax\in\mathbb{F}_{q^n}[x]$ is a permutation
polynomial if and only if $(-1)^{n/d}a^{{(q^n-1)}/{(q^{d}-1)}}\neq
1$ where $d=(n,r)$. In this paper, the compositional inverse of
$L_r(x)$ is explicitly determined when this condition holds.
\end{abstract}


\end{frontmatter}

\section{Introduction}\label{secintro}
Let $\mathbb{F}_{q}$ be the finite field with $q$ elements where $q$
is a prime or a prime power. A polynomial over $\mathbb{F}_{q}$ is
called a permutation polynomial  if it can induce  a bijective map
from $\mathbb{F}_{q}$ to itself. For a given permutation polynomial
$f(x)$ over $\mathbb{F}_{q}$, the polynomial $f^{-1}(x)$ that can
induce the inverse map of the map induced by $f(x)$ is called the
compositional inverse of $f(x)$. That is to say, $f^{-1}(x)$
satisfies the relation
\[f(f^{-1}(x))\equiv f^{-1}(f(x))\equiv x\mod
(x^q-x).\] Generally speaking, it is far form a simple matter to
explicitly represent compositional inverses of known classes of
permutation polynomials over finite fields. Up to present, there are
only a few classes of permutation polynomials whose compositional
inverses can be determined. We refer to \cite{lidlDic,coulter,wu2},
for example, for some results on this topic.

Linearized polynomials are of special interest in studying
permutation polynomials over finite fields. A linearized polynomial
over the finite field $\mathbb{F}_{q^n}$ is a polynomial of the form
\[L(x)=\sum_{i=0}^{n-1}a_ix^{q^i}.\]
A well-known criterion of Dickson says that $L(x)$ is a permutation
polynomial over $\mathbb{F}_{q^n}$ if and only if $\det D_L\neq0$
\cite{lidl}, where
\[D_L=\begin{pmatrix}
a_0&a_1&\dots&a_{n-1}\\
a_{n-1}^q&a_0^q&\dots&a_{n-2}^q\\
\vdots&\vdots&&\vdots\\
a_1^{q^{n-1}}&a_2^{q^{n-1}}&\dots&a_0^{q^{n-1}}
\end{pmatrix}\]
is called the associate Dickson matrix of $L(x)$ \cite{wu}. As a
supplement of this result, in \cite{wu} the author and Liu found a
relation between the compositional inverse of a linearized
permutation polynomial and the inverse of its associate Dickson
matrix, obtaining the following result.

\begin{thm}[See \cite{wu}]\label{invLPP}
Let $L(x)\in\mathbb{F}_{q^n}[x]$ be a linearized permutation
polynomial. Then $D_{L^{-1}}=D_L^{-1}$. More precisely, if
$L(x)=\sum_{i=0}^{n-1}a_ix^{q^i}$ and $\tilde{a}_i$ is the
$(i,0)$-th cofactor of $D_L$, $0\leq i\leq n-1$, then
\[L^{-1}(x)=\frac{1}{\det L}\sum_{i=0}^{n-1}\tilde{a}_ix^{q^i},\]
where $\det
L=\sum_{i=0}^{n-1}a_{n-i}^{q^i}\tilde{a}_i\in\mathbb{F}_{q}$
(subscripts reduced modulo $n$).
\end{thm}

Though Dickson's criterion is simple, it is not convenient enough to
use to characterize linearized permutation polynomials sometimes.
This is because singularity of certain matrices over finite fields
cannot be easily characterized, especially when some entries of them
contain parameters. Also, Theorem \ref{invLPP} only presents a
general method to determine compositional inverses of linearized
permutation polynomials and in some cases it cannot be conveniently
utilized. The main difficulty lies in computing determinants of
certain matrices over finite fields. Therefore, compositional
inverses of linearized permutation polynomials is not easy to obtain
in general. For example, in \cite{wu3} the author explicitly
determined the compositional inverses of a class of linearized
permutation polynomials of  simple forms over $\mathbb{F}_{2^n}$ for
an odd $n$, but the  computations were rather complicated.

The simplest examples of linearized polynomials over
$\mathbb{F}_{q^n}$ are linearized monomials, which are always
permutation ones unless it is zero. Moreover, the compositional
inverse of a nonzero linearized monomial can be trivially
determined. Non-trivial examples of linearized polynomials that
should be considered firstly are linearized binomials of the form
\[L_{r,s}(x)=ax^{q^r}+bx^{q^s},~a,~b\in\mathbb{F}_{q^n},~1\leq s<r\leq n-1.\]
However, when studying permutation behavior of such linearized
polynomials, it is easy to see that the problem can be reduced to
studying linearized binomials of the form
\begin{equation}\label{LBP}
L_r(x)=x^{q^r}+ax,~a\in\mathbb{F}_{q^n},~1\leq r\leq n-1.
\end{equation}
It is well known that the condition under which $L_r(x)$ in
\eqref{LBP} is a linearized permutation polynomial over
$\mathbb{F}_{q^n}$ is
\begin{equation}\label{condition}
(-1)^{n/d}a^{{(q^n-1)}/{(q^{d}-1)}}\neq 1
\end{equation}  where $d=(n,r)$.
Actually, $L_r(x)$ can permute $\mathbb{F}_{q^n}$ if and only if $0$
is the only solution of $L_r(x)=0$ in $\mathbb{F}_{q^n}$, or
equivalently, $x^{q^r-1}=-a$ has no solution in $\mathbb{F}_{q^n}$.
This is further equivalent to
\[-a\not\in(\mathbb{F}_{q^n}^*)^{q^r-1}=
\mathcal {C}_{\frac{q^n-1}{(q^{n}-1,q^r-1)}}=\mathcal
{C}_{\frac{q^n-1}{q^{d}-1}},\] i.e.
\[(-a)^{\frac{q^n-1}{q^{d}-1}}=(-1)^{\frac{n}{d}}a^{\frac{q^n-1}{q^{d}-1}}\neq 1,\]
where $\mathcal{C}_l$ denotes the subgroup of $\mathbb{F}_{q^n}^*$
of order $l$ for any $l$ with $l\mid (q^n-1)$. However, noticing
that the associate Dickson matrix of $L_r(x)$ is
\[D_{L_r}=\begin{pmatrix}
a&&&&&1&&&\\
&a^q&&&&&1\\
&&\ddots&&&&&\ddots\\
&&&a^{q^{n-1-r}}&&&&&1\\
1&&&&a^{q^{n-r}}\\
&1&&&&a^{q^{n-r+1}}\\
&&1&&&&a^{q^{n-r+2}}\\
&&&\ddots&&&&\ddots\\
&&&&1&&&&a^{q^{n-1}}
\end{pmatrix},\]
whose shape  is flexible since $r$ is not fixed, we find it
difficult to compute $\det D_{L_r}$ to derive condition
\eqref{condition} for $L_r(x)$ to be a permutation polynomial by
Dickson's criterion. Furthermore, it is more difficult to compute
cofactors of elements in the first column of $D_{L_r}$ to obtain
$L_r^{-1}(x)$ by Theorem \ref{invLPP} when condition
\eqref{condition} holds.

In the sequel we will always assume condition \eqref{condition}
holds and  devote to getting explicit representation of the
compositional inverse of $L_r(x)$ for any $1\leq r\leq n-1$. Our
main obversion is that the problem can be reduced to getting
compositional inverse of $L_1(x)$, which can be easily solved from a
direct utilization of Theorem \ref{invLPP}. The idea of the
reduction process may be applied in other similar problems, which
will be generally discussed in the end.

\section{Compositional inverses of linearized permutation
binomials}\label{secinvLBP}

Denote by ``$\nor_{n:l}$" the norm map from $\mathbb{F}_{q^n}$ to
$\mathbb{F}_{q^l}$ for some $l\mid n$, i.e.
$\nor_{n:l}(x)=x^{(q^n-1)/(q^l-1)}$ for any $x\in\mathbb{F}_{q^n}$,
and when $l=1$,  denote by ``$\nor$" the absolute norm map for
simplicity. Then condition \eqref{condition} can be rewritten as
$(-1)^{n/d}\nor_{n:d}(a)\neq1$. The following theorem gives explicit
representation of $L_r^{-1}(x)$.

\begin{thm}\label{invLBP}
\[L_r^{-1}(x)=\frac{\nor_{n:d}(a)}{\nor_{n:d}(a)+(-1)^{\frac{n}{d}-1}}\sum_{i=0}^{\frac{n}{d}-1}(-1)^i
a^{-\frac{q^{(i+1)r}-1}{q^r-1}}x^{q^{ir}}.\]
\end{thm}
\begin{pf}
It is obvious that $\nor_{n:d}(a)^{q^r}=\nor_{n:d}(a)$ as $d\mid r$.
Besides, since $(n/d,r/d)=1$, we have
\[1+q^d+\ldots+q^{(\frac{n}{d}-1)d}\equiv 1+q^r+\ldots+q^{(\frac{n}{d}-1)r}\mod (q^n-1),\]
and thus
\[a^{-\frac{q^{\frac{n}{d}r}-1}{q^r-1}}=\nor_{n:d}\left(\frac{1}{a}\right).\]
For any $x\in\mathbb{F}_{q^n}$, we have
\begin{eqnarray*}
  &&L_r(L_r^{-1}(x))\\ &=&
\frac{\nor_{n:d}(a)}{\nor_{n:d}(a)+(-1)^{\frac{n}{d}-1}}\left[
  \sum_{i=0}^{\frac{n}{d}-1}(-1)^i
a^{-\frac{q^{(i+1)r}-1}{q^r-1}q^r}x^{q^{(i+1)r}}+\sum_{i=0}^{\frac{n}{d}-1}(-1)^i
a^{-\frac{q^{(i+1)r}-1}{q^r-1}+1}x^{q^{ir}}\right]\\
  &=&\frac{\nor_{n:d}(a)}{\nor_{n:d}(a)+(-1)^{\frac{n}{d}-1}}\left[
  \sum_{i=0}^{\frac{n}{d}-1}(-1)^i
a^{-\frac{q^{(i+2)r}-q^r}{q^r-1}}x^{q^{(i+1)r}}+\sum_{i=0}^{\frac{n}{d}-1}(-1)^i
a^{-\frac{q^{(i+1)r}-q^r}{q^r-1}}x^{q^{ir}}\right] \\
   &=& \frac{\nor_{n:d}(a)}{\nor_{n:d}(a)+(-1)^{\frac{n}{d}-1}}\left[(-1)^{\frac{n}{d}-1}\nor_{n:d}\left(\frac{1}{a}\right)^{q^r}
   x^{\frac{n}{d}r}+\sum_{i=1}^{\frac{n}{d}-1}(-1)^{i-1}
a^{-\frac{q^{(i+1)r}-q^r}{q^r-1}}x^{q^{ir}}\right.\\
&&\qquad\qquad\qquad\qquad~~~~\left.+(-1)^0x+\sum_{i=1}^{\frac{n}{d}-1}(-1)^i
a^{-\frac{q^{(i+1)r}-q^r}{q^r-1}}x^{q^{ir}}\right] \\
   &=&\frac{\nor_{n:d}(a)}{\nor_{n:d}(a)+(-1)^{\frac{n}{d}-1}}\left[\frac{(-1)^{\frac{n}{d}-1}}{\nor_{n:d}(a)}+1\right]x  \\
   &=&x.
\end{eqnarray*}
\qedd
\end{pf}

For some special cases of $r$, we can directly get the following
corollaries.

\begin{cor}\label{corr=1}
$L_1(x)$ is a linearized permutation polynomial over
$\mathbb{F}_{q^n}$ if and only if $(-1)^n\nor(a)\neq 1$, and under
this condition, we have
\[L_1^{-1}(x)=\frac{\nor(a)}{\nor(a)+(-1)^{n-1}}\sum_{i=0}^{n-1}(-1)^i
a^{-\frac{q^{i+1}-1}{q-1}}x^{q^{i}}.\]
\end{cor}

\begin{cor}\label{corcoprime}
Assume $(r,n)=1$. Then $L_{r}(x)$ is a linearized permutation
polynomial over $\mathbb{F}_{q^n}$ if and only if $(-1)^n\nor(a)\neq
1$, and under this condition, we have
\[L_r^{-1}(x)=\frac{\nor(a)}{\nor(a)+(-1)^{n-1}}\sum_{i=0}^{n-1}(-1)^i
a^{-\frac{q^{(i+1)r}-1}{q^r-1}}x^{q^{ir}}.\]
\end{cor}

\begin{cor}
Assume $n$ is even. Then $L_{n/2}(x)$ is a linearized permutation
polynomial over $\mathbb{F}_{q^n}$ if and only if $a^{q^{n/2}+1}\neq
1$, and under this condition, we have
\[L_{\frac{n}{2}}^{-1}(x)=\frac{1}{a^{q^{\frac{n}{2}}+1}-1}
\left(a^{q^{\frac{n}{2}}}x-{x}^{q^{\frac{n}{2}}}\right).\]
\end{cor}

\section{The method to obtain Theorem
\ref{invLBP}}\label{sectechnuque}
In this section, we explain the detail of our method to derive
$L_r^{-1}(x)$. As mentioned in Section \ref{secintro}, we cannot use
Theorem \ref{invLPP} directly since the shape of $D_{L_r} $ is
flexible. Now, by some trikes,  we reduce the problem to one that
can be easily handled.

Firstly, we let $q_1=q^d$ and then $L_r(x)=x^{q_1^{r/d}}+ax$, which
can be viewed as a linearized polynomial over
$\mathbb{F}_{q_1^{n/d}}=\mathbb{F}_{q^n}$. Since $(r/d,n/d)=1$, this
implies that we need only to consider the initial problem in the
case $(r,n)=1$;

Secondly, when $(r,n)=1$, we let $q_2=q^r$. Consider the composite
field of $\mathbb{F}_{q^n}$ and $\mathbb{F}_{q_2}$, which is just
$\mathbb{F}_{q^{nr}}=\mathbb{F}_{q_2^{n}}$, and view
$L_r(x)=x^{q_2}+ax$ as a linearized polynomial over
$\mathbb{F}_{q_2^n}$. Since $(r,n)=1$, we have
\[ 1+q_2+\ldots+q_2^{n-1}\equiv1+q+\ldots+q^{n-1}\mod (q^n-1),
\]
hence when $(-1)^na^{{(q^n-1)}/{(q-1)}}\neq1$, we know that
\[(-a)^{\frac{q_2^n-1}{q_2-1}}=(-a)^{1+q^r+\ldots+q^{(n-1)r}}=(-1)^na^{1+q+\ldots+q^{n-1}}=
(-1)^na^{\frac{q^n-1}{q-1}}\neq 1,\] i.e.
\[a\not\in\mathcal{C}_{\frac{q_2^n-1}{q_2-1}}=\left(\mathbb{F}_{q_2^{n}}^*\right)^{q_2-1}.\]
This implies that $L_r(x)$ can induce a permutation of
$\mathbb{F}_{q_2^{n}}$. Furthermore, since
$L_r(\mathbb{F}_{q^n})\subseteq\mathbb{F}_{q^n}$ and $L_r(x)$ can
induce a permutation of $\mathbb{F}_{q^{n}}$, which is a subset of
$\mathbb{F}_{q_2^{n}}$, the compositional inverse of $L_r(x)$ viewed
as a linearized permutation polynomial over $\mathbb{F}_{q_2^{n}}$
must be the the compositional inverse of $L_r(x)$ viewed as a
linearized permutation polynomial over $\mathbb{F}_{q^{n}}$, after
reduction modulo $(x^{q^n}-x)$. To this end, we  need only to
consider the initial problem in the case $r=1$.

To summarize, if we can determine $L_1^{-1}(x)$, i.e. if we can
obtain the result of Corollary \ref{corr=1} at first, then we can
determine $L_r^{-1}(x)$ for $(r,n)=1$  via replacing $q$ by
$q_2=q^r$ in the representation of $L_1^{-1}(x)$, obtaining the
result of Corollary \ref{corcoprime} (note that the relation
$\nor_{nr:r}(a)=\nor(a)$ should be used). Afterwards, for a general
$1\leq r\leq n-1$ with $(r,n)=d$, we can replace $q$, $r$ and $n$ by
   $q_1=q^d$, $r/d$ and $n/d$, respectively, in the representation
  of $L_r^{-1}(x)$ in Corollary \ref{corcoprime}, to obtain the
  representation of $L_r^{-1}(x)$ in Theorem \ref{invLBP}.

The only rest problem is to derive $L_1^{-1}(x)$. This can be done
via directly using Theorem \ref{invLPP} since in this case $D_{L_1}$
has a fixed shape, namely,
\[D_{L_1}=\begin{pmatrix}
a&1\\
&a^q&1\\
&&\ddots&\ddots\\
&&&a^{q^{n-2}}&1\\
1&&&&a^{q^{n-1}}
\end{pmatrix}.\]
Note that the $(0,0)$-th cofactor of it is
\[\tilde{a}_0=(-1)^{1+1}\det\begin{pmatrix}
a^q&1\\
&\ddots&\ddots\\
&&a^{q^{n-2}}&1\\
&&&a^{q^{n-1}}
\end{pmatrix}=\frac{\nor(a)}{a},\]
the $(i,0)$-th cofactor of it is
\begin{eqnarray*}
  \tilde{a}_i &=&(-1)^{i+1+1}\det \begin{pmatrix}
1\\
a^q&1\\
&\ddots&\ddots\\
&&a^{q^{i-1}}&1\\
&&&&a^{q^{i+1}}&1\\
&&&&&\ddots&\ddots\\
&&&&&&a^{q^{n-2}}&1\\
&&&&&&&a^{q^{n-1}}
\end{pmatrix} \\
   &=& (-1)^i\frac{\nor(a)}{a^{1+q+\ldots+q^i}} \\
   &=&(-1)^i\nor(a)a^{-\frac{q^{i+1}-1}{q-1}}
\end{eqnarray*}
for any $1\leq i\leq n-2$, and the $(n-1,0)$-th cofactor of it is
\[\tilde{a}_{n-1}=(-1)^{n+1}\det\begin{pmatrix}
1\\
a^q&1\\
&\ddots&\ddots\\
&&a^{q^{n-2}}&1
\end{pmatrix}=(-1)^{n-1}.\]
Finally we have
\[\det L_1=a\tilde{a}_0+\tilde{a}_{n-1}=\nor(a)+(-1)^{n-1}.\]
Then Corollary \ref{corr=1} follows.

\section{A general discussion}\label{secgeneraldiss}
In fact, the technique we introduc in Section \ref{sectechnuque} to
reduce the problem to one that is simple enough to deal with is
enlightening. Other  problems related to linearized polynomials over
$\mathbb{F}_{q^{n}}$, especially those whose terms are all of the
form $x^{q^{it}}$, can be similarly handled. Note that the main
observation we make in the second step of the problem reduction
process is that, $L_r(x)$ can induce a permutation of
$\mathbb{F}_{q^{rn}}$ if it can induce a permutation of
$\mathbb{F}_{q^{n}}$. This fact can be affirmatively generalized.

\begin{thm}\label{generalfact}
Let $t$ and $n$ be positive integers with $(t,n)=1$, and
$\bar{q}=q^t$ for a prime power $q$. Assume the linearized
polynomial $L(x)=\sum_{i=0}^{n-1}a_ix^{q^i}\in\mathbb{F}_{q^n}[x]$
can induce a permutation of $\mathbb{F}_{q^{n}}$. Then the
linearized polynomial
$\bar{L}(x)=\sum_{i=0}^{n-1}a_{ti}x^{\bar{q}^i}$ (subscripts reduced
modulo $n$) can induce a permutation of $\mathbb{F}_{\bar{q}^n}$.
\end{thm}

\begin{rmk}
From $(t,n)=1$ we know that the composite field of
$\mathbb{F}_{q^n}$ and $\mathbb{F}_{\bar{q}}$ is
$\mathbb{F}_{q^{nt}}=\mathbb{F}_{\bar{q}^{n}}$. Besides, it is
obvious that $\bar{L}(\mathbb{F}_{q^n})\subseteq\mathbb{F}_{q^n}$,
and\[\bar{L}(x)=\sum_{i=0}^{n-1}a_{ti}x^{q^{ti}}\equiv L(x)\mod
(x^{q^n}-x),\]thus $\bar{L}(x)$ can permute $\mathbb{F}_{q^n}$ as
well. Therefore, Theorem \ref{generalfact} actually presents a
method to extend a linearized permutation polynomial over a ``small"
field  to be a linearized permutation polynomial over a ``big"
field. See the graph of field extensions below.
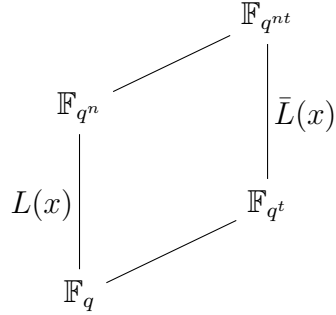
\begin{center}
\begin{tikzpicture}
\path (0,0) node(a){$\mathbb{F}_{q}$} (0,2.5)
node(b){$\mathbb{F}_{q^n}$} (2.5,1.2) node(c){$\mathbb{F}_{q^t}$}
(2.5,3.7) node(d){$\mathbb{F}_{q^{nt}}$}; \draw
(a)--(b)--(d)--(c)--(a); \path (-.5,1.2)
node[draw=none,fill=none](e){$L(x)$} (3,2.4)
node[draw=none,fill=none](f){$\bar{L}(x)$}; \draw (e) (f);
\end{tikzpicture}
\end{center}
\end{rmk}

To prove Theorem \ref{generalfact}, the following lemmas are needed.

\begin{lem}[See \cite{mene}]\label{lembasis}
Let $t$ and $n$ be positive integers with $(t,n)=1$, and
$\{\beta_i\}_{i=0}^{n-1}$ be a basis of $\mathbb{F}_{q^n}$ over
$\mathbb{F}_{q}$. Then $\{\beta_i\}_{i=0}^{n-1}$ is a basis of
$\mathbb{F}_{q^{nt}}$ over $\mathbb{F}_{q^t}$.
\end{lem}

\begin{lem}[See \cite{wu}]\label{lemlpprep}
Let $L(x)=\sum_{i=0}^{n-1}a_ix^{q^i}$ be a linearized permutation
polynomial over $\mathbb{F}_{q^n}$. Then there exist two bases
$\{\alpha_i\}_{i=0}^{n-1}$ and $\{\beta_i\}_{i=0}^{n-1}$ of
$\mathbb{F}_{q^n}$ over $\mathbb{F}_{q}$ such that
\[(a_0,a_1,\ldots,a_{n-1})=(\alpha_0,\alpha_2,\ldots,\alpha_{n-1})\left(\beta_i^{q^j}\right)_{0\leq i,j\leq n-1}.\]
\end{lem}

\noindent\textbf{Proof of Theorem \ref{generalfact}.} Since $L(x)$
is a linearized permutation polynomial over $\mathbb{F}_{q^n}$,
there exist two bases $\{\alpha_i\}_{i=0}^{n-1}$ and
$\{\beta_i\}_{i=0}^{n-1}$ of $\mathbb{F}_{q^n}$ over
$\mathbb{F}_{q}$ such that
\[(a_0,a_1,\ldots,a_{n-1})=(\alpha_0,\alpha_2,\ldots,\alpha_{n-1})\left(\beta_i^{q^j}\right)_{0\leq i,j\leq n-1}\]
according to Lemma \ref{lemlpprep}. As $(t,n)=1$, we have
\[(a_0,a_t,\ldots,a_{(n-1)t})=(\alpha_0,\alpha_2,\ldots,\alpha_{n-1})\left(\beta_i^{q^j}\right)_{0\leq i,j\leq n-1}P,\]
where $P$ is a permutation matrix with entries 1 in the $(it,i)$-th
place for $0\leq i\leq n-1$ and 0 in other places. Note that
\[\left(\beta_i^{q^j}\right)_{0\leq i,j\leq n-1}P=\left(\beta_i^{\bar{q}^j}\right)_{0\leq i,j\leq n-1}.\]
From Lemma \ref{lembasis} we know  $\{\alpha_i\}_{i=0}^{n-1}$ and
$\{\beta_i\}_{i=0}^{n-1}$  are two bases of $\mathbb{F}_{\bar{q}^n}$
over $\mathbb{F}_{\bar{q}}$, thus we finally get that $\bar{L}(x)$
is a linearized  permutation polynomial over
$\mathbb{F}_{\bar{q}^n}$ applying Lemma \ref{lemlpprep} again.\qedd

\section{Concluding remarks}

In this paper,  the explicit representation of a linearized
permutation binomial of the form $L_r(x)=x^{q^r}+ax$ over the finite
fields $\mathbb{F}_{q^n}$ is derived. Our main tool is Theorem
\ref{invLPP}, but it can only be used after we find by some trikes
that talking about the initial problem for $L_1(x)$ is enough. We
should point out that,  though it cannot be conveniently utilized
sometimes, Theorem \ref{invLPP} is quite useful in studying
linearized permutation polynomials of special types and their
compositional inverses. In fact, we have also used it to compute
compositional inverses of certain linearized permutation trinomials
over finite fields, the results of which will be proposed in a
further paper.



\begin{thebibliography}{115}



\bibitem{charpin}
P. Charpin, G. Kyureghyan, When does $G(x)+\gamma
\mbox{\rm{Tr}}(H(x))$ permutate $\mathbb{F}_{p^n}$?, Finite Fields
Appl. 15 (2009) 615--632.

\bibitem{coulter}
R.S. Coulter, M. Henderson, The compositional inverse of a class of
permutation polynomials over a finite field, Bull. Austral. Math.
Soc. 65 (2002) 521--526.


\bibitem{lidl}
R. Lidl, H. Niederreiter, Finite fields, second edn., Encyclopedia
Math. Appl., vol. 20, Cambridge University Press, New York, 1997.

\bibitem{lidlDic}
R. Lidl, G.L. Mullen, G. Turnwald, Dickson polynomials, Pitman
monographs and surveys in pure and applied mathematics, vol. 65,
Longman Scientific \& Technical, Essex,  1993.

\bibitem{mene}
A. Menezes, I. Blake, X. Gao, et al., Applications of Finite Fields,
Kluwer Academic, Boston, 1993.



\bibitem{wu}
B. Wu, Z. Liu, Linearized polynomials over finite fields revisited,
Finite Fields Appl. 22 (2013) 79--100.

\bibitem{wu2}
B. Wu, Z. Liu, The compositional inverse of a class of bilinear
permutation polynomials over finite fields of characteristic 2,
Finite Fields Appl 24 (2013) 136--147.

\bibitem{wu3}
B. Wu, The compositional inverse of a class of linearized
permutation polynomials over $\mathbb{F}_{2^n}$, $n$ odd,
arXiv:1305.1411v2 [math.CO], preprint (2013).

\end{thebibliography}
\end{document}